\documentclass{amsproc}

\usepackage{amssymb}

\usepackage{graphicx}

\usepackage{}

\newtheorem{theorem}{Theorem}[section]

\theoremstyle{definition}

\theoremstyle{remark}

\numberwithin{equation}{section}

\newtheorem{problem}[theorem]{Problem}
\newtheorem{construction}[theorem]{Construction}
\newtheorem{conjecture}[theorem]{Conjecture}

\newcommand{\dmrjdel}[1]{}

\newcommand{\Aut}{\mathrm{Aut}\,}

\newcommand{\pd}{\partial}
\newcommand{\Jnorm}[2]{\left\langle J_{#1}, J_{#1} \right\rangle_{#2} }

\newcommand{\Pbar}{\mathsf{Pic}}

\newcommand{\cC}{\mathcal{C}}

\newcommand{\sF}{\mathsf{F}}
\newcommand{\bi}{\mathbf{i}}
\newcommand{\sK}{\mathsf{K}}

\newcommand{\cM}{\mathcal{M}}
\newcommand{\cP}{\mathcal{P}}

\newcommand{\bbR}{\mathbb{R}}
\newcommand{\frS}{\mathfrak{S}}

\newcommand{\bx}{\mathbf{x}}
\newcommand{\by}{\mathbf{y}}
\newcommand{\bz}{\mathbf{z}}

\newcommand{\al}{\alpha}
\newcommand{\be}{\beta}

\newcommand{\la}{\lambda}
\newcommand{\La}{\Lambda}
\newcommand{\vep}{\varepsilon}

\begin{document}

\title{Transitive factorizations of permutations and geometry}


\author{I.P.Goulden}
\address{Dept. of Combinatorics and Optimization, University of Waterloo, Canada}
\curraddr{}
\email{ipgoulde@uwaterloo.ca}
\thanks{The work of both authors was supported by NSERC Discovery Grants.}

\author{D.M.Jackson}
\address{Dept. of Combinatorics and Optimization, University of Waterloo, Canada}
\curraddr{}
\email{dmjackso@uwaterloo.ca}


\date{}

\begin{abstract} 
We give an account of our work on transitive factorizations of permutations. The work has had impact upon other areas of mathematics such as the enumeration of graph embeddings, random matrices, branched covers, and the moduli spaces of curves.  Aspects of these seemingly unrelated areas are seen to be related in a unifying view from the perspective of algebraic combinatorics. At several points this work has intertwined with Richard Stanley's in significant ways.
\end{abstract}

\maketitle


\begin{center}
\textit{Dedicated to Richard Stanley on the occasion of his 70th birthday.}
\end{center}

\vspace{.1in}

\section{Introduction} \label{Intro}

Richard Stanley has been a pioneer in modern combinatorics, and a key figure in the development of both enumerative combinatorics and algebraic combinatorics.  

In enumerative combinatorics two crucial building blocks are \textbf{(a)}  the generating series for a set of combinatorial objects, and \textbf{(b)}~the relationship between algebraic operations on types of generating series and combinatorial operations on the set.  The enumerative significance of a generating series is, of course, that the normalized coefficient of each of its monomials counts the objects in the set of combinatorial objects indexed by the monomial.  Stanley contributed early to these important building blocks in his paper with Doubilet and Rota~\cite{drs} -- part of Gian-Carlo Rota's seminal series \emph{``On the foundations of combinatorial theory''}. Further early work appeared in paper~\cite{s1}.  

The power of algebraic combinatorics often seems to depend on the efficacy of the analogue relationship between algebra and the combinatorics, in which methods from one may assist in solving questions raised in the other.  Stanley has been particularly attracted by combinatorics that has made an impact in other branches of mathematics, and was himself an early developer of many of these analogue relationships.

Our own work has been inspired by these early developments.  In an essential way, they have influenced our work in enumerative combinatorics, both together and separately, which has, in its turn,  contributed to the further study of the connection between combinatorial structure and algebraic structure, and its application to other parts of mathematics and the mathematical sciences.  

\section{Transitive factorizations of permutations} 

In this article, we describe our longtime work on transitive factorizations of permutations. The themes that it illustrates include:
\begin{itemize}
\item the fundamental underlying combinatorial problem is very simple to state;
\item the contexts in which instances of this combinatorial problem arise are diverse within mathematics and mathematical physics;
\item the interplay between algebra and combinatorics is exhibited in both directions, with methods from other parts of mathematics applied to combinatorial problems, as well as combinatorial methods applied within other parts of mathematics;
\item work in this area continues to be the subject of intense research activity both in algebraic combinatorics and in other parts of mathematics;
\item Stanley's work has made an important contribution in a number of places.
\end{itemize}

We now describe the fundamental factorization problem that we consider in this article, with two variations. The following notation will be used:  $\frS_n$ is the symmetric group acting on the symbols $\{1,\ldots,n\}$; we write $\al\vdash n$, or equivalently $|\al |=n$, to indicate that $\al$ is a partition of $n$; the number of parts in $\al$ is denoted by $l(\al )$; $\cC_\al$ is used to denote the conjugacy class of $\frS_n$ with natural index $\al$. If $m_i$ is the number of parts of $\al$ equal to $i$, $i\geq 1$, then $|\Aut \al |=\prod_{i\geq 1} m_i!$.

\begin{problem}[{\bf The Permutation Factorization Problem}]\label{pfp}
For fixed partitions $\al, \be_1, \ldots, \be_m$ of $n$, find the number of  permutations $\rho\in\cC_{\al}$ and $\pi_i\in\cC_{\be_i}$, $i=1,\ldots ,m$, such that
\begin{equation}\label{permfactn}
\pi_1 \pi_2 \cdots \pi_m  = \rho.
\end{equation}
We shall call $(\pi_1,\ldots, \pi_m)$ a \emph{factorization} of $\rho$.
\end{problem}

\begin{problem}[{\bf The Transitive Permutation Factorization Problem}]\label{tpfp}
For fixed partitions $\al, \be_1, \ldots, \be_m$ of $n$, find the number of  permutations $\rho\in\cC_{\al}$ and $\pi_i\in\cC_{\be_i}$, $i=1,\ldots ,m$, that satisfy equation~(\ref{permfactn}), and such that $\langle\pi_1,\ldots ,\pi_m\rangle$, the group generated by the factors $\pi_1,\ldots, \pi_m$, acts transitively on the underlying symbols $\{ 1,\ldots ,n\}$.
In this case we shall call $(\pi_1,\ldots, \pi_m)$ a \emph{transitive factorization} of $\rho$.
\end{problem}

There were a number of papers in the combinatorics literature on permutation factorization problems in the 70's by various authors. These relied on elementary methods only; see, for example, Walkup~\cite{wa}. In his 1981 paper~\cite{s2}, Stanley applied the powerful mathematical methodology of symmetric group characters to solve the problem in the case in which all factors are $n$-cycles (in the conjugacy class $\cC_{(n)}$). As part of this, he was able to prove a conjecture from~\cite{wa}.

A convenient way of describing the method of symmetric group characters is to work in the centre of the group algebra of $\frS_{n}$. One basis of the centre is the set 
$\{\sK_{\theta} = \sum_{\sigma\in\cC_{\theta}} \sigma\; \colon \theta\vdash n\}$ of classes, and another is the set $\{\sF_\al \colon \al\vdash n\}$ of orthogonal idempotents. These bases are related by the linear relations
\begin{equation}\label{classidem}
\sF_\al = \frac{\chi^{\al}(1^n)}{n!} \sum_{\theta\vdash n} \chi^{\al}(\theta) \sK_{\theta},\qquad\qquad \sK_{\theta} = |\cC_{\theta}| \sum_{\al\vdash n} \frac{\chi^\al(\theta)}{\chi^{\al}(1^n)} \sF_{\al},
\end{equation}
where  $\chi^\al(\theta)$ is the character $\chi^\al$ of the (ordinary) irreducible representation of $\frS_{n}$ indexed by $\al$, and evaluated on the class $\cC_\theta$.

Encoded in this way, the answer to the Permutation Factorization Problem is given by
\begin{equation}\label{classsoln}
|\cC_{\al}|\cdot\Big( [\sK_{\al}]\sK_{\be_1}\cdots\sK_{\be_m}\Big) ,
\end{equation}
where the notation $[X]Y$ denotes the coefficient of $X$ in the expansion of $Y$. The factor $|\cC_{\al}|$ appears in~(\ref{classsoln}) since each element of the class $\cC_{\al}$ is created in the product with the same frequency; the factor $|\cC_{\al}|$ would be removed if in~(\ref{permfactn}) we were considering permutation factorizations of a fixed and arbitrary element $\rho$ of the class $\cC_{\al}$. Of course, to apply~(\ref{classsoln}), one simply applies~(\ref{classidem}), and uses the fact that $\sF_{\al}\cdot\sF_{\be} = \sF_{\al}$ if $\al =\be$, and $\sF_{\al}\cdot\sF_{\be} = 0$ otherwise. Thus, one has changed bases to one in which multiplication is ``trivial'', before changing back to the basis of conjugacy classes. In general, the resulting expression is a sum over partitions of $n$ involving arbitrary characters of the symmetric group. Such summations are generally regarded as intractable, but significant simplification occurs in the case considered by Stanley~\cite{s2}, where all factors are $n$-cycles (so $\be_i=(n)$ for $i=1,\ldots ,m$).
In this case, the characters have explicit evaluations, almost always equal to $0$.

Since the group generated by any single $n$-cycle acts transitively on $\{1,\ldots,n\}$, the factorizations that Stanley considered in that paper were in fact transitive, though this condition is not caused by any particularly ``natural'' mathematical reason.  The remainder of the paper deals with applications of transitive permutation factorization in which the transitivity condition is quite natural, and involve factors in arbitrary conjugacy classes, not simply $n$-cycles.

\section{Maps in orientable surfaces}\label{maporient} 

\begin{figure}
\scalebox{.4}{\includegraphics{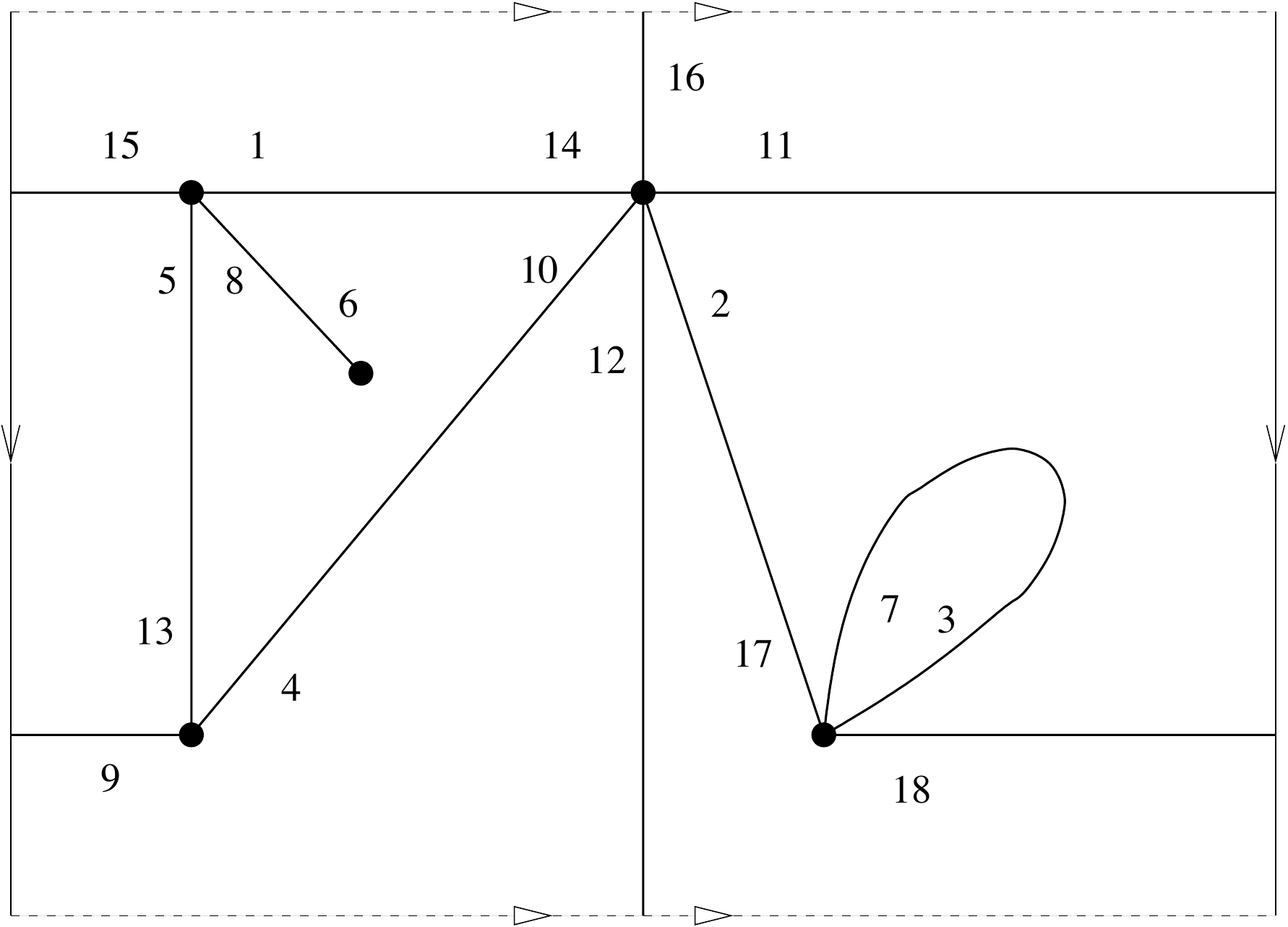}}
\caption{A decorated rooted map embedded in the torus}\label{otblefig}
\end{figure}

A \emph{rooted map} is a graph embedded in a surface so that all faces are two-cells (homeomorphic to a disc). In the case of orientable surfaces, one vertex is distinguished, called the root vertex, and one edge incident with the root vertex is distinguished, called the root edge. In order to construct permutation factorizations corresponding to a rooted map in an orientable surface with $n$ edges, assign labels to the two ends of the edges with the integers $1,\ldots ,2n$ subject only to the restriction that the end of the root edge incident with the root vertex is assigned the label $1$. We call the resulting object a \emph{decorated rooted map}, and of course, there are $(2n-1)!$ decorated rooted maps corresponding to every rooted map with $n$ edges. An example with $9$ edges using the standard polygonal representation of the torus is given in Figure~\ref{otblefig}. 
  
 \begin{construction} \label{mapperms}
(see, \textit{e.g.}, Tutte~\cite{t} for full details) Given a decorated rooted map with $n$ edges, construct three permutations $\nu ,\vep ,\phi $ in $\frS_{2n}$ as follows:

\begin{itemize}
\item
the disjoint cycles of $\nu$, the \emph{vertex permutation}, are the clockwise circular lists of end labels of edges incident with each vertex;

\item
the disjoint cycles of  $\vep$, the \emph{edge permutation}, are the pairs of labels on the two ends of each edge;

\item
the disjoint cycles of $\phi$, the \emph{face permutation}, are the counterclockwise circular lists of the second label on each edge encountered when traversing the interior of the faces. 
 
\end{itemize}
\end{construction}

As an example of Construction~\ref{mapperms}, the three permutations that we construct from the decorated rooted map given in Figure~\ref{otblefig} are:
\begin{align*}
\nu &= (1\, 8\, 5\, 15)(1\, 12\, 10\, 14\, 16\, 11)(3\, 18\, 17\, 7)(4\, 9\, 13)(6),\\
\epsilon &= (1\, 14)(2\, 17)(3\, 7)(4\, 10)(5\, 13)(6\, 8)(9\, 18)(11\, 15)(12\, 16),\\
\phi &= (1\, 6\, 8\, 13\, 10)(2\, 16\, 15\, 14\, 12\, 4\, 18)(3\, 9\, 5\, 11\, 17)(7).
\end{align*}

From the description in Construction~\ref{mapperms}, it is clear that in general, as in the above example,
\begin{itemize}

\item
the lengths of the cycles of $\nu$ specify the vertex-degrees of the underlying rooted map,

\item
all cycles of $\vep$ have length $2$,

\item
the lengths of the cycles of $\phi$ specify the face-degrees of the underlying map.

\end{itemize}

\noindent
Moreover, by construction we have $\vep\nu = \phi$, and the fact that $\langle \vep , \nu\rangle$ acts transitively on the symbols $1,\ldots ,2n$ follows immediately from the fact that the embedded graph is connected. Finally, the genus of the embedding surface can be obtained from $\nu ,\vep ,\phi$ by Euler's formula.

Consequently, the enumeration of rooted maps embedded in orientable surfaces is, up to scaling, a special case of the Transitive Permutation Factorization Problem (Problem~\ref{tpfp}), in which there are precisely two factors. When we solve this enumerative question in terms of group characters by means of~(\ref{classsoln}), and form the generating series, we find that symmetric functions are introduced in a natural way because the linear relations~(\ref{classidem}) are scale equivalent to the linear relations
 \begin{equation*}
s_{\al}= \sum_{\theta\vdash n}\frac{|\cC_{\theta}|}{n!} \chi^{\al}(\theta) p_{\theta}, \qquad\qquad p_{\theta} =\sum_{\al\vdash n} \chi^{\al}(\theta) s_{\al},
\end{equation*}
between the Schur functions $s_{\al}$ and power sums $p_{\theta}$. Thus, for $\la$, $\mu$ partitions of $2n$, if $m^{\la}_{\mu}$ is the number of rooted maps in orientable surfaces with $n$ edges, vertex degrees given by the parts of $\la$, and face degrees given by the parts of $\mu$, then we obtain
\begin{equation*}\label{coefforient}
m^{\la}_{\mu} = [p_{\la}(\bx)p_{\mu}(\by)p_{(2^n)}(\bz)t^{2n}] H_O\left( p(\bx), p(\by), p(\bz), t \right) ,
\end{equation*}
where
\begin{equation}\label{gseriesorient}
H_O\left( p(\bx), p(\by), p(\bz), t \right) =  
t\frac{\pd}{\pd t} \log \left( \sum_{\theta\in\cP} 
\frac{|\theta|!}{\chi^\theta (1^{|\theta|})} \,
s_\theta(\bx) s_\theta(\by) s_\theta(\bz) \, t^{|\theta|} \right),
\end{equation}
and $\cP$ is the set of all (integer) partitions, $p(\bx) := (p_1(\bx), p_2(\bx), \ldots)$,  $p_k(\bx)$ is the degree $k$ power sum symmetric function in the indeterminates $\bx =(x_1, x_2, \ldots)$. Full details of this was developed with Visentin in \cite{jv1, jv2}, so we make only a few technical remarks here: \textbf{(a)} the generating series $H_O$ is actually an \emph{exponential} generating series in the indeterminate $t$, but for the number of \emph{decorated} rooted maps, \textbf{(b}) the ``$\log$'' appears in~(\ref{gseriesorient}) to restrict to the connected objects in the usual way for exponential generating series,  \textbf{(c}) the effect of $t\pd/\pd t$ in~(\ref{gseriesorient}) is to multiply the coefficient of $t^{2n}$ by $2n$, thus adjusting the exponential monomial $\frac{t^{2n}}{(2n)!}$ to $\frac{t^{2n}}{(2n-1)!}$; this division by $(2n-1)!$ is the  correct scaling between decorated rooted maps and rooted maps, \textbf{(d)} the coefficient of arbitrary monomials $p_{\tau}(\bz )$ in $H_O$ also has combinatorial meaning; it accounts for rooted \emph{hypermaps}.

\section{Maps in surfaces and Jack symmetric functions}

This enumerative approach to rooted maps was extended in~\cite{gj4} from orientable surfaces to all surfaces (includes non-orientable surfaces). For all surfaces, the class algebra of the symmetric group -- products of conjugacy classes, was replaced by the Hecke algebra associated with the hyperoctahdedral group -- products of double cosets of the symmetric group multiplied by the hyperoctahedral subgroup on both sides.  Stanley's paper with Hanlon and Stembridge~\cite{hss} was an essential source, describing completely the character theory of this algebra, and the relationship with symmetric functions, in this case the zonal polynomials $Z_{\theta}$.
For the generating series, again with $\la$, $\mu$ partitions of $2n$, let $\ell^{\la}_{\mu}$ be the number of rooted maps in locally orientable surfaces with $n$ edges, vertex degrees given by the parts of $\la$, and face degrees given by the parts of $\mu$. Then we obtain
\begin{equation*}\label{coefflocorient}
{\ell}^{\la}_{\mu} = [p_{\la}(\bx)p_{\mu}(\by)p_{(2^n)}(\bz)t^{2n}] H\left( p(\bx), p(\by), p(\bz), t \right) ,
\end{equation*}
where
\begin{equation*}\label{gserieslocorient}
H\left( p(\bx), p(\by), p(\bz), t \right) =  
2t\frac{\pd}{\pd t} \log \left( \sum_{\theta\in\cP} 
\frac{\chi^{2\theta}(1^{|2\theta|})}{|2\theta|!} \,
Z_\theta(\bx) Z_\theta(\by) Z_\theta(\bz) \, t^{|\theta|} \right) ,
\end{equation*}
and $2\theta:=(2\theta_1, 2\theta_2,\ldots)$ for $\theta=(\theta_1, \theta_2, \ldots)$. Again, the coefficient of arbitrary monomials $p_{\tau}(\bz )$ in $H$ accounts for rooted \emph{hypermaps}.

But there is more that we can say.  We showed in~\cite{gj3}  that $H_O$ and $H$ have a common generalization as the cases $\al =1$ and $\al=2$, respectively, of
\begin{equation}\label{Jacksum}
\Psi(\bx, \by, \bz; t,\al )
:= \al  \, t \frac{\pd}{\pd t}  \log  \sum_{\theta\in\cP} \frac{1}{\Jnorm{\theta}{\al }}
J_\theta(\bx; \al ) \, J_\theta(\by; \al ) \, J_\theta(\bz; \al )\, t^{|\theta|} ,
\end{equation}
where $J_\theta(\bx; \al )$ is the Jack symmetric function with parameter $\al$ and $\langle\cdot,\cdot\rangle_{\al}$ is the standard inner product on the ring of symmetric functions. 
In this work, our source for the necessary results on Jack symmetric functions was again a paper of Stanley, in this case~\cite{s3}. 
Following extensive computer algebra computations with the generating series $\Psi$, we conjectured  the following: 
\begin{conjecture}[{\bf The $b$-conjecture}]\label{bconjecture}
The series $\Psi(\bx, \by, \bz; t,1+b)$ has coefficients that are \emph{polynomial} in $b$ with \emph{non-negative integer coefficients}. In this polynomial, the constant term, obtained with $b=0$ (so $\al =1$), accounts for rooted hypermaps embedded in orientable surfaces, and the sum of all terms, obtained with $b=1$ (so $\al =2$), accounts for rooted hypermaps embedded in all surfaces. Accordingly, the indeterminate $b$ marks a \emph{statistic of nonorientability} associated with rooted hypermaps.
\end{conjecture}

The $b$-conjecture has not yet been resolved, but some progress towards determining a suitable statistic of nonorientabilty has been made, and work is ongoing. Recent progress with La Croix~\cite{JlaC} has involved providing combinatorial interpretations, in terms of maps and hypermaps (or of transformations of the series in terms of polynomial glueings), for sums of coefficients rather than for individual coefficients. Despite success with these marginal sums, a complete understanding of the $b$-Conjecture for maps and hypermaps continues to elude us.
 
There is a closely related $b$-conjecture for \emph{matchings} that we conjectured in \cite{gj3}. Previously, progress on the matching version has appeared in~Do{l}ega and F\'{e}ray~\cite{df} and Do{l}ega, F\'{e}ray and \'{S}niady~\cite{dfs}.

\section{Maps, matrix integrals and virtual Euler characteristic}

Stanley's paper with Hanlon and Stembridge~\cite{hss} also contained matrix integral results associated with symmetric functions. Especially due to the influence of mathematical physics, such results are important in algebraic combinatorics, and this continues to be an area of active research interest. The generating series for maps, whose symmetric function expressions have been discussed in the previous two sections, also have matrix integral forms, and we present a brief discussion of these in this section. 

For $\bi=(i_1,i_2,\ldots )$, let $m_O(\bi,j,n)$ denote the number of rooted maps in orientable surfaces with $i_k$ vertices of degree $k$, $k\ge 1$, $j$ faces, and $n$ edges, and let $m(\bi,j,n)$ denote the corresponding number in all surfaces. Define the generating series
$$
M_O(\by,x,z)  = \sum_{\bi, j, n} m_O(\bi,j,n) \, \by^\bi x^j z^n
\quad\mbox{and}\quad
M(\by,x,z) = \sum_{\bi, j, n} m(\bi,j,n) \, \by^\bi x^j z^n
$$
where $\by^\bi := \prod_{k\ge1} y_k^{i_k}$. A matrix integral, over Hermitian complex matrices, was given in \cite{j2} for $M_O$. Using a different argument, a matrix integral over real symmetric matrices was given
in \cite{gj65} for $M$. These integrals can be transformed by the Weyl integration theorems
(the diagonalizing groups are the unitary group and the orthogonal group, respectively; the measure may be factored into Haar measure for the manifold of the groups, and an $\bbR^N$ integral over the spectra $\la$), and with Harer~\cite{ghj} we obtained a common generalization $M(\by,N,z; \al )$ of the diagonalized integrals, where $M(\by,x,z; 1 )= M_O(\by,x,z)$ and $M(\by,x,z; 2 )= M(\by,x,z)$, with
\begin{equation}\label{RNalphaintegral}
M(\by,N,z; \al )  
=2\al z\frac{\pd}{\pd z}\log 
\left(
\frac
{\int_{\bbR^N} \left|V(\la)\right|^{\frac{2}{\al}}e^{\sum_{k\ge1} \frac{1}{k} y_k \sqrt{z}^k p_k(\la)} 
\cdot e^{- \frac{1}{2\al}  p_2(\la)}  d\la}
{\int_{\bbR^N} \left|V(\la)\right|^{\frac{2}{\al}} e^{- \frac{1}{2\al} p_2(\la)}  d\la}
\right) ,
\end{equation}
and $V$ is the Vandermonde determinant. For combinatorial reasons, the coefficients of $z^n$ are polynomials in $N$; we may formally replace $N$ by $x$ to obtain $M(\by,x,z;\al )$ from $M(\by,N,z;\al )$.

Note that the parameter $\al$ in~(\ref{RNalphaintegral}) specializes in the same way as the Jack parameter in~(\ref{Jacksum}) or Conjecture~\ref{bconjecture}, but we do not have a matrix integral (undiagonalized) that involves the parameter $\al$.

An immediate application in~\cite{ghj} was to obtain $\chi(\cM^s_g(\tau))$, the virtual Euler characteristic of the moduli spaces of real curves of genus $g$, with $s$ marked points and a fixed topological type of orientation reversing involution $\tau$.  Harer and Zagier~\cite{hz} had earlier obtained $\chi(\cM^s_g)$, the virtual Euler characteristic for the case of complex curves, using the fact that it can be obtained from a sum over rooted \emph{monopoles} -- maps with a single face. A further application in \cite{ghj}, of the common generalization~(\ref{RNalphaintegral}), was to determine a common generalization $\xi^s_g(\al )$ of these virtual Euler characteristics, that gave the complex case when $\al =1$ and the real case when $\al =2$. Comparing this with the $b$-conjecture (Conjecture~\ref{bconjecture}) suggests that the coefficients of $b$ in the polynomial $\xi^s_g(1+b)$ have a geometric interpretation in the context of the moduli spaces of curves, but this has not been resolved to date.

There are applications of $\xi^s_g(\al )$ to string theory. For example, the expressions for the virtual Euler characteristics for real and complex curve confirm the case $g=1$ determined by Ooguri and Vafa~\cite{ov} associated with the $SO(N)$ and $Sp(N)$ gauge groups.

We conclude this section with the following observation. By comparing the symmetric function and matrix integral expressions for the map generating series, we conjectured in \cite{gj65} that
$$
\left\langle J_\theta(\la;\al)\right\rangle_{\bbR^N} = J_\theta(1_N;\al) \cdot [p_2^m]\,J_\theta, \mbox{where}\; \left\langle f(\la)\right\rangle_{\bbR^N} := \frac {\int_{\bbR^N} |V(\la)|^{\frac{2}{\al}} e^{-\frac{1}{2\al} p_2(\la)} f(\la) d\la} {\int_{\bbR^N} |V(\la)|^{\frac{2}{\al}}  e^{-\frac{1}{2\al} p_2(\la)}  d\la},
$$
and $\theta\vdash 2m$. This conjecture was subsequently proved by Okounkov~\cite{ok}.

\section{Branched covers of the sphere and Hurwitz numbers} \label{S:EncTranFact}

In Section~\ref{maporient}, we showed that the special case of the Transitive Factorization Problem (Problem~\ref{tpfp}) with two factors has a geometric interpretation in terms of rooted maps, or embedded graphs, in orientable surfaces. In that case the genus of the embedding surface could be determined from the factors by Euler's polyhedral formula. In this Section, we consider a second geometric interpretation of the Transitive Factorization Problem, in this case in terms of \emph{branched covers} from algebraic geometry. 

Consider branched covers of the sphere by an $n$-sheeted Riemann surface of genus $g$. Suppose that the branch points are $P_0,P_1,\ldots, P_m$, with branching at $P_i$ specified by permutation $\pi_i\in\frS_n$, for $i=0,1,\ldots ,m$, where $\pi_0\in\cC_{\al}$ and $\pi_i\in\cC_{\be_i}$, $i=1,\ldots ,m$. (This means that if one walks in a small neighbourhood, counterclockwise, around $P_i$, starting at sheet $j$, then one ends at sheet $\pi_i(j)$.) Hurwitz~\cite{h} proved that, up to homeomorphism, each $\pi_0,\pi_1,\ldots ,\pi_m$ as defined above determines a unique branched cover precisely when, in the language of Problem~\ref{tpfp}, $(\pi_1,\ldots ,\pi_m)$ is a transitive factorization of $\rho =\pi_0^{-1}$. Note the following points:

\begin{itemize}

\item
the fact that the permutations form a factorization is a \emph{monodromy} condition on the sheets;

\item
the transitivity condition on the factorization means that the cover is connected;

\item
we say that the \emph{branching type} of $P_0$ is $\al$, and of $P_i$ is $\be_i$, $i=1,\ldots ,m$;

\item
the genus of the surface $g$ is obtained from the branching types of the permutations by the \emph{Riemann-Hurwitz} formula, which gives
\begin{equation}\label{RieHur}
\sum_{i=1}^m \left( n-l(\be_i)\right) = n+l(\al)+2g-2;
\end{equation}

\item
the minimum number of factors in such a transitive factorization, from~(\ref{RieHur}), is $n+l(\al)-2$ which are obtained with genus $g=0$. We call such factorizations \emph{minimal} transitive factorizations;

\item 
if branching at a branch point is a transposition then it is called \emph{simple}.

\end{itemize}

Our own work on the enumeration of branched covers was initiated through Richard Stanley. Arising from joint work with Crescimanno \cite{ct}, Washington Taylor (Dept. of Physics, MIT) had asked Stanley about a particular transitive factorization problem for permutations, that turned out to be a special case of \emph{Hurwitz numbers} in genus $0$. Stanley suggested he should contact me (DMJ). Taylor's e-mail languished unanswered for three months on an old main frame computer at Waterloo.  It was only through Stanley's well-known and encyclop{\ae}dic grasp of progress on active questions that I became aware of the oversight, after he e-mailed asking about progress.

The Hurwitz number $H_{\al}^g$ is the number of topologically distinct branched covers in genus $g$, in which branching is of type $\al$ at one specified branch point, and branching is simple at $r$ remaining branch points. \emph{Topologically distinct} means that we divide the number of branched covers by $n!$, for geometric reasons. Thus $H_{\al}^g$ equals $\frac{1}{n!}$ times the number of transitive factorizations of an element of $\cC_{\al}$ into $r$ transpositions, where:
\begin{itemize}

\item
from~(\ref{RieHur}), the number of transpositions is given by $r=n+l(\al)+2g-2$;

\item
the group generated by the $m$ transpositions acts transitively on $\{ 1,\ldots ,n\}$. Equivalently, the multigraph with vertex-set $\{ 1,\ldots ,n\}$, and $r$ edges, one edge $\{ a,b\}$ for each transposition $(a\, b)$, is connected.
\end{itemize}

In this language, Taylor was asking about the number of transitive factorizations of the identity permutation into $2n-2$ transpositions. These are minimal transitive factorizations, with genus $g=0$, and hence are given by the Hurwitz number $H_{(1^n)}^0$. 

\section{The join-cut equation}

In our first attempts to solve Taylor's problem, we applied group characters, to obtain a generating series in the form of a logarithm of a Schur function summation, analogous to the map generating series given in~(\ref{gseriesorient}). We were not able to obtain an explicit formula for Taylor's problem from this symmetric function form of the generating series, so we moved on to the following more indirect analysis: Form the generating series
\begin{equation*}
H^0=\sum_{n=1}^{\infty}z^n \sum_{\al\vdash n} \frac{H_{\al}^0}{(n+l(\al )-2)!} p_{\al}
\end{equation*}
in the indeterminates $z,p_1,p_2,\ldots$, and suppose that the last transposition in the factorization is $(a\, b)$. Then when we multiply a permutation by $(a\, b)$, there are two possibilities:  
\begin{itemize}

\item[\textbf{Case 1:}]  $a$ and $b$ occur on \emph{different} cycles of lengths $i$ and $j$, and the cycles are \emph{joined} to form a single cycle of length $i+j$;

\item[\textbf{Case 2:}]  $a$ and $b$ occur on the \emph{same} cycle, of length $i+j$, and this cycle is \emph{cut} into a cycle of length $i$ and a cycle of length $j$.

\end{itemize}
This analysis (which we call a join-cut analysis)  leads immediately to the formal partial differential equation
\begin{equation}\label{jcutpde}
\frac{1}{2}\sum_{i,j\ge 1} \left( p_{i+j} \,  i\frac{\pd H^0}{\pd p_i} j\frac{\pd H^0}{\pd p_j} +  p_i p_j  (i+j)\frac{\pd H^0}{\pd p_{i+j}} \right)
-z\frac{\pd H^0}{\pd z} - \sum_{i\ge 1} p_i \frac{\pd H^0}{\pd p_i} + 2H^0 = 0,
\end{equation}
which we call the \emph{join-cut equation} for the series $H^0$. Together with the initial condition $[z^0]H^0=0$, this uniquely determines $H^0$.

It turns out that working with equation~(\ref{jcutpde}) is greatly simplified by changing variables from $z$ to $s$ by means of the functional equation
\begin{equation}\label{szeqn}
s=z\exp\left( \sum_{i\ge 1}\frac{i^i}{i!}p_i\, s^i  \right) .
\end{equation}
This change of variables is pefectly natural within algebraic combinatorics, and is quite tractable by means of Lagrange's Implicit Function Theorem. For example, we can express the series $H^0$ in terms of $s$ in the simple form
\begin{equation*}
\left( \frac{\pd^2}{\pd z^2}\right) H^0 = \log \left( \frac{s}{z} \right) ,
\end{equation*}
and it follows immediately from Lagrange's Theorem that the Hurwitz number for genus $0$ has the explicit form
\begin{equation}\label{H0alpha}
H^0_{\al}=  \frac{(n+\ell -2)!}{|\Aut\al |} \, n^{\ell -3} \, \prod_{j=1}^{\ell } \frac{\al_j^{\al_j}}{\al_j!},
\end{equation}
where $\al\vdash n$ and $\ell = l(\al )$ (see \cite{gj1} for full details).  In this notation,  D\'{e}nes~\cite{d1} and Crescimanno and Taylor~\cite{ct} had previously obtained the results for $\al = (n)$ and $\al=(1^n)$, respectively. We were unaware when writing the paper that this explicit form for all $\al$ had been obtained much earlier by Hurwitz~\cite{h}.

The join-cut analysis can be extended to Hurwitz numbers in arbitrary genus. Again the change of variables in~(\ref{szeqn}) helps to simplify, and in~\cite{gjvn} (see also \cite{gj7} and \cite{gjv1}) we were led to conjecture  the existence of a polynomial $P_{g,\ell}$, for each $g\ge 0$ and $\ell\geq 1$, such that for all partitions $\alpha \vdash n$ with $\ell=l(\al )$ parts,
\begin{equation}\label{Hgalpha} 
  H^g_{\al} = \frac{(n+\ell +2g-2)!}{|\Aut\al |} \, P_{g,\ell}(\alpha_1, \dots, \alpha_\ell) \, \prod_{j=1}^\ell \frac{\alpha_j^{\alpha_j}}{\al_j!}.
\end{equation}
For example, from~(\ref{H0alpha}), since $\al\vdash n$, we have
\begin{equation*}
P_{0,\ell}(\al_1,\ldots ,\al_{\ell}) = \left( \al_1+\cdots +\al_{\ell} \right)^{\ell -3}.
\end{equation*}

\section{Hodge integrals, the moduli space of curves, and integrable hierarchies}

Hurwitz numbers have been the subject of much research interest over the last couple of decades, with a variety of mathematical areas making substantial contributions, including mathematical physics, algebraic geometry and algebraic combinatorics. For example, soon after we conjectured the existence of the polynomial $P_{g,\ell }$ in~(\ref{Hgalpha}), Ekedahl, Lando, Shapiro and Vainshtein \cite{elsv} proved it by constructing an explicit expression for the polynomial as a \emph{Hodge integral}. The expression is the celebrated ELSV formula
\begin{equation}\label{eq:elsvf}
  P_{g,\ell}(\alpha_1, \dots, \alpha_\ell) = \int_{{\overline{\cM}}_{g,\ell}} \frac{1 - \lambda_1 + \cdots + (-1)^g \lambda_g}{(1 - \alpha_1 \psi_1) \cdots (1 - \alpha_\ell \psi_\ell)},
\end{equation}
where ${\overline{\cM}}_{g,\ell}$ is the (compact) moduli space of stable $\ell$-pointed genus $g$ curves, $\psi_1$, $\dots$, $\psi_{\ell}$ are (codimension $1$) classes corresponding to the $\ell$ marked points, and $\lambda_k$ is the (codimension $k$) $k$th Chern class of the Hodge bundle. Equation~(\ref{eq:elsvf}) should be interpreted as follows: formally invert the denominator of the integrand as a geometric series; select the terms of codimension $\dim {\overline{\cM}}_{g,\ell}=3g - 3 + \ell$; and ``intersect'' these terms on ${\overline{\cM}}_{g,\ell}$. 

Earlier, and perhaps most notably, Witten~\cite{wi} had initiated much of this work by his conjecture that a transform of the generating series for Hurwitz numbers is a $\tau$-function for the KdV hierarchy from integrable systems. Witten's motivation for the conjecture was that two different models of two-dimensional quantum gravity have the same partition function. For one of these models, the partition function can be described in terms of intersection numbers on moduli space, but also in terms of Hurwitz numbers. Witten's Conjecture \cite{wi} was proved soon after by Kontsevich \cite{kon}, and a number of proofs have appeared since, for example Kazarian and Lando \cite{kl}.

A variant of Hurwitz numbers called \emph{double} Hurwitz numbers have also been the subject of recent research interest. They were introduced by Okounkov \cite{ok2}, motivated by a conjecture of Pandharipande \cite{p} in Gromow-Witten theory. Okounkov showed that a particular generating series for double Hurwitz numbers is a $\tau$-function for the Toda lattice hierarchy from integrable systems.

The double Hurwitz number $H_{\al ,\be}^g$ is the number of topologically distinct branched covers in genus $g$, where branching is of type $\al$ at one specified branch point, type $\be$ at another specified branch point, and branching is simple at $r$ remaining branch points. Thus $H_{\al ,\be}^g$ equals $|\Aut\al |\cdot |\Aut \be |/n!$ (this factor is chosen for geometric reasons) times the number of transitive factorizations of an element of $\cC_{\al}$ into an element of $\cC_{\be}$ together with $r$ transpositions, where:
\begin{itemize}

\item
from~(\ref{RieHur}), the number of transpositions is given by $r=l(\al)+l(\be)+2g-2$;

\item
the group generated by the element of $\cC_{\be}$ and the $r$ transpositions acts transitively on $\{ 1,\ldots ,n\}$.
\end{itemize}

In joint work with Vakil \cite{gjv2}, we used both group characters and a join-cut analysis to obtain various results for double Hurwitz numbers. One of these was the following conjectured ELSV-type formula for double Hurwitz numbers where one of the partitions has a single part:

$$
H^g_{\al ,(n)} = n\, (\ell +2g-1)! \,\int_{\overline{\Pbar}_{g,\ell }}
\frac{\La_0 -\La_2 + \cdots \pm\La_{2g}}
{(1-\al_1\psi_1) \cdots (1-\al_{\ell}\psi_{\ell})} ,
$$
where $\al =(\al_1,\ldots ,\al_{\ell})$, $\overline{\Pbar}_{g,n}$ is a conjectural compactification of the universal Picard variety, and $\La_{2k}$ is a conjectural (codimension $2k$) class.

A second result we obtained for double Hurwitz numbers, reminiscent of the polynomiality result given in~(\ref{Hgalpha}) for Hurwitz numbers, was a \emph{piecewise} polynomiality result. In particular, for fixed $g,\ell ,k$, and $\al =(\al_1,\ldots ,\al_{\ell})$ with $\ell$ parts, and $\be =(\be_1,\ldots ,\be_k)$ with $k$ parts, then $H^g_{\al ,\be}$ is piecewise polynomial (and \emph{not} polynomial) in the parts $\al_1,\ldots ,\al_{\ell},\be_1,\ldots ,\be_k$, of degree $4g-3+\ell -k$. Our proof of Piecewise Polynomiality used ribbon graphs to interpret double Hurwitz numbers as counting lattice points in certain polytopes. We then required Ehrhart's Theorem and Ehrhart polynomials, whose properties have been studied extensively by Stanley, see for example \cite{s4}. The piecewise polynomiality property of the double Hurwitz numbers has prompted further study of the chamber structure and wall crossings in these polytopes; see for example \cite{cjm, ssv}. 

 Finally, a substantially different but related geometric setting in which transitive permutation factorizations have been applied is given in Lando and Zvonkin \cite{lz}, where they are called \emph{constellations}.


\bibliographystyle{amsplain}

\end{document}